\documentclass[11pt,a4paper,lualatex]{amsart}

\usepackage[marginparwidth=0pt,margin=20truemm]{geometry} 

\usepackage{mypackage} 
\usepackage{mycommand} 
\usetikzlibrary{knots}
\usetikzlibrary{positioning}

\newcommand{\fg}{\mathfrak{g}}
\newcommand{\su}{\mathfrak{su}}
\renewcommand{\iu}{i}

\usepackage[pdfencoding=auto,hypertexnames=false]{hyperref}
\hypersetup{colorlinks=false}

\numberwithin{equation}{section} 
\usepackage{mytheoremeng} 

\begin{document}

\title[Homological blocks with simple Lie algebras and WRT invariants]{Homological blocks with simple Lie algebras and Witten--Reshetikhin--Turaev invariants}

\author[Y. Murakami]{Yuya Murakami}
\address{Faculty of Mathematics, Kyushu University, 744, Motooka, Nishi-ku, Fukuoka, 819-0395, JAPAN}
\email{murakami.yuya.896@m.kyushu-u.ac.jp}

\author[Y. Terashima]{Yuji Terashima}
\address{Graduate school of science, Tohoku University, 
	6-3, Aoba, Aramaki-aza, Aoba-ku, 
	Sendai, 980-8578, Japan}
\email{yujiterashima@tohoku.ac.jp}

\date{\today}

\maketitle


\begin{abstract}
	In this article, for any Seifert fibered homology $3$-sphere, we introduce homological blocks with simple Lie algebra and prove that its radial limits are identified with the Witten--Reshetikhin--Turaev invariants.
	To prove it, we develop an asymptotic formula and a vanishing result of asymptotic coefficients.
\end{abstract}


\tableofcontents


\section{Introduction} \label{sec:intro}


Homological blocks are $q$-series introduced by Gukov--Pei--Putrov--Vafa~\cite{GPPV} for plumbed $3$-manifolds based on \cite{GPV} 
and are very interesting objects in both physics and mathematics.
Their properties has been much studied recently~\cite{CCFGH,Chun,Chung_rational,Chung_resurgent,GHNPPS,GM,GMP,GPV,Wu}.
One of the most important properties that are expected from a physical viewpoint is that their special limits at roots of unity 
are identified with the Witten--Reshetikhin--Turaev (WRT) invariants. This is a mathematical conjecture and is solved in~\cite{FIMT} and~\cite{Andersen-Mistegard} for Seifert fibered homology $3$-spheres and in~\cite{M_GPPV} for negative definite plumbed manifolds generalizing ~\cite{MM,M_plumbed}. See also \cite{LZ, LR,  H_Bries,H_Lattice,H_Seifert,H_Lattice2} for important related works done before homological blocks were introduced.

The original homological blocks are associated with the Lie algebra $\su(2)$. It is hoped that there is a definition of homological blocks with general Lie algebras that have several good properties like the original homological blocks. In this article, introducing homological blocks for Seifert fibered homology $3$-spheres with simple Lie algebras $\fg$, we prove that their special limits at roots of unity are identified with the WRT invariants when the Lie algebra $\fg$ is simply-laced. This identification is proved by comparing our asymptotic formulas of generalized homological blocks which is a main result in this paper to formulas of WRT invariants with simply laced simple Lie algebra $\fg$ by Mari\~{n}o \cite{Marino} which generalize formulas by Lawrence--Rozansky \cite{LR} for $\fg=\su(2)$. See \cite{Hansen-Takata_asymptotic,Hansen-Takata} for earlier works on WRT invariants of Seifert fibered homology $3$-spheres with simple Lie algebras. Remark that our generalized homological blocks coincide with the original homological blocks when the Lie algebra $\fg$ is $\su(2)$. Also, we find that our generalized homological blocks for $\fg=\su(N)$ coincide with the generalization defined in \cite{Chung_SU(N)}. 
Remark that Park~\cite{Park,EGGKPS}  gives a generalization for negative definite plumbed 3-manifolds.
We hope that our definition agrees with his definition for Seifert fibered homology spheres, and can generalize our results to negative definite plumbed 3-manifolds.

Our proof of the asymptotic fomulas of generalized homological blocks consists of two steps.
In the first step, we develop a new asymptotic formula of infinite series whose idea is based on \cite[p.~98, Proposition]{LZ}, \cite[Equation (44)]{Zagier_asymptotic}, \cite[Equation (2.8)]{BKM}, \cite[Lemma 2.2]{BMM_high_depth} and \cite[Proposition 5.6 and Lemma 5.11]{M_plumbed}.
In the final step, we prove a vanishing result of asymptotic coefficients by considering holomorphy of a certain meromorphic function.
We deal with it by combining ideas of \cite[Lemma 6]{FIMT} and \cite[Proposition 6.1]{M_plumbed}.

It is expected that homological blocks have good modular transformation properties. In fact, the paper \cite{Matsusaka-Terashima} obtained explicit modular transformation formulas of original homological blocks for any Seifert fibered homology $3$-sphere, and uses them to give a new proof of the Witten asymptotic conjecture.  
It would be interesting to get generalizations for simple Lie algebras and prove the Witten asymptotic conjecture for simple Lie algebras that still has not been verified. 

This article is organized as follows: 
In \cref{sec:homological_block}, we prepare the settings for Seifert fibered homology spheres, define their generalized homological blocks and state our main result.
In \cref{sec:asymp}, we develop asymptotic formulas corresponding to the two steps mentioned above.
In \cref{sec:proof_ex}, we prove that the WRT invariants are identified with radial limits of our generalized homological blocks.
Moreover, we give examples of our homological blocks, and consequently, it becomes clear that our generalized homological blocks coincide with homological blocks defined in \cite{Chung_SU(N)} for $\fg = \su(N)$.
In \cref{sec:WRT}, as a consequence of our main result, we prove that WRT invariants can be expressed as radial limits of our homological blocks.


\section*{Acknowledgement}
We would like to thank H. Fuji, K. Hikami, K. Iwaki, M. Kaneko, T. Matsusaka, W. Misteg\aa{}rd, A. Mori, H. Murakami, S. Sugimoto and T. Yamauchi for valuable discussions. This work is partially supported by
JSPS KAKENHI Grant Number JP23KJ1675, 
JP21K03240 
and JP22H01117. 


\section{Homological blocks with simple Lie algebras} \label{sec:homological_block}


In this section, we define homological blocks with simple Lie algebras and states our main result.

To begin with, we prepare settings for Seifert fibered homology spheres.

\begin{figure}[tb]
	\centering
	\begin{tikzpicture}
		\begin{knot}[
			flip crossing/.list={2,4,6}
			]
			\strand (1,0) circle[x radius=6cm, y radius=2cm];
			\strand[] (-2,2) circle[x radius=0.5cm, y radius=1.5cm];
			\strand[] (0,2) circle[x radius=0.5cm, y radius=1.5cm];
			\strand[] (4,2) circle[x radius=0.5cm, y radius=1.5cm];
			
			\node (1) at (1, -2.5) {$ 0 $};
			\node (2) at (-2, 4) {$ p_1 / q_1 $};
			\node (3) at (0, 4) {$ p_2 / q_2 $};
			\node (4) at (1.5, 4) {$ \cdots $};
			\node (5) at (4, 4) {$ p_n / q_n $};
		\end{knot}
	\end{tikzpicture}
	\caption{The surgery diagram of $ M(p_1/q_1, \dots, p_n/q_n) $.} \label{fig:surgery_diagram}
\end{figure}
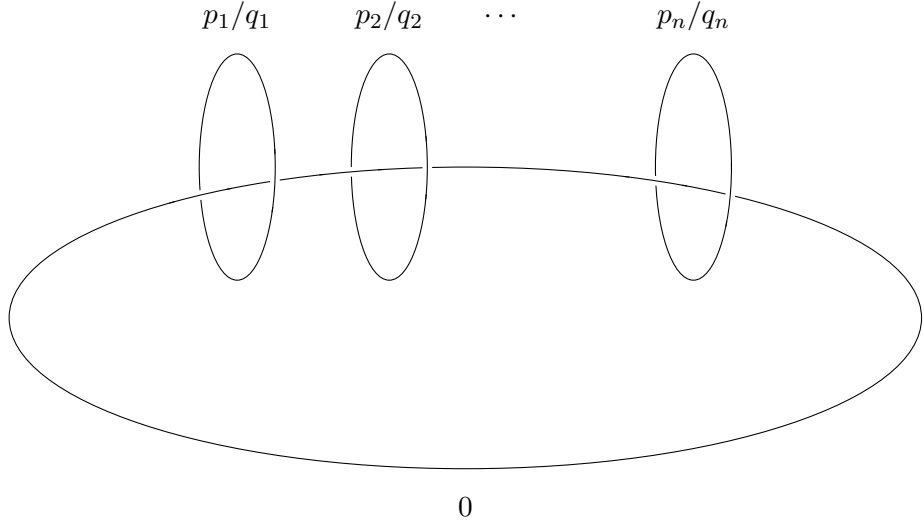

Let $ n \ge 3 $ be an integer and $ (p_{1}, q_{1}), \dots, (p_{n}, q_{n}) $ be pairs of coprime integers such that $ 0 < q_1 < p_1, \dots, 0 < q_n < p_n $, $ p_1, \dots, p_n $ are pairwise coprime, and
\[
P \sum_{1 \le i \le n} \frac{q_i}{p_i} = 1,
\]
where $ P := p_1 \cdots p_n $.
We denote by $ M := M(p_1/q_1, \dots, p_n/q_n) $ the Seifert fibered $ 3 $-manifold with $ n $-singular fibers and surgery integers $ p_1, \dots, p_n $ and $ q_1, \dots, q_n $.
This manifold is obtained by the surgery diagram in \cref{fig:surgery_diagram}.

Let $ \frakg $ be a complex simple Lie algebra and $ \Delta_{+} $ be a set of positive roots.
Then, we define the homological block for $ M $ with $ \frakg $ as
\[
\Phi(q)
:=
q^{-(\dim_{\bbC} \frakg) \phi \abs{\rho}^{2}/2}
\sum_{m \in \Z_{\ge m_0}^{\Delta_+}} q^{\abs{\sum_{\alpha \in \Delta_+} m_\alpha \alpha }^{2} / 8P}
\prod_{\alpha \in \Delta_+} \chi_{m_\alpha},
\]
where 
\begin{itemize}
	\item denote
	\[
	\phi := 3 - \frac{1}{P} + 12 \sum_{i=1}^n s(q_i, p_i),
	\]
	where 
	\[
	s(q_i, p_i) := \frac{1}{4q_i} \sum_{j=1}^{q_i - 1}
	\cot \left( \frac{\pi j}{q_i} \right) \cot \left( \frac{\pi j p_i}{q_i} \right)
	\]
	are the Dedekind sums,
	\item denote $ \rho := \sum_{\alpha \in \Delta_{+}} \alpha/2 \in \frakh_\R^* $, 
	\item define $ \abs{ \cdot } $ as follows:
	\begin{itemize}
		\item let $ \frakh $ be the Cartan subalgebra of $ \frakg $;
		\item let $ \frakh_\R^* $ be the real vector subspace of the dual space $ \frakh^* $ generated by the root system $ \Delta $;
		\item let $ \sprod{\cdot, \cdot} $ be the standard inner product on $ \frakh_\R^* $ normalized such that the longest roots have length $ \sqrt{2} $;
		\item for $ x \in \frakh_\R^* $, denote $ \abs{x} := \sqrt{\sprod{x, x}} $;
	\end{itemize}	
	\item denote
	\[
	\sum_{m = m_0}^\infty \chi_m q^{m/2P}
	:=
	G_{p_1, \dots, p_n} (q)
	:=
	(q^{1/2} - q^{-1/2})^{2-n} \prod_{1 \le i \le n} (q_\alpha^{1/2p_i} - q_\alpha^{-1/2p_i}).
	\]
\end{itemize}

\begin{rem} \label{rem:chi_m}
	We can express $m_0$ and $\chi_m$ as
	\begin{align}
		m_0 &= P \left( - \frac{1}{p_1} - \cdots - \frac{1}{p_n} + n - 2 \right), \\
		\chi_m &=
		\begin{dcases}
			(-1)^n \veps_1 \cdots \veps_n \binom{m'+n-3}{n-3},
			& 
			\begin{gathered}
				\text{ if } \frac{m}{P} = \frac{\veps_1}{p_1} + \cdots + \frac{\veps_n}{p_n} + n - 2 + 2m'
				\text{ for some } \\
				\veps_1, \dots, \veps_n \in \{ \pm 1 \} \text{ and } m' \in \Z_{\ge 0},
			\end{gathered}
			\\
			0 & \text{ otherwise}.
		\end{dcases}
	\end{align}
\end{rem}

\begin{rem}
	The homological block $\Phi(q)$ is independent of a choice of a set of positive roots $ \Delta_{+} $.
	This is because, if we choose an other set of positive roots $ \Delta'_{+} $, then we have an element $w$ of the Weyl group $W$ such that $ \Delta_{+} = w \Delta'_{+} $ and the standard inner product $ \sprod{\cdot, \cdot} $ is invariant under the action of the Weyl group $W$.
\end{rem}

Our main theorem is the following statement.

\begin{thm} \label{thm:main}
	For each positive integer $ k $, we have
	\begin{align}
		\lim_{q \to \zeta_k} \Phi(q)
		= \,
		&\bm{e} \left( 
		-\frac{(\dim_{\bbC} \frakg) \phi \abs{\rho}^{2} }{2k}
		\right)
		\left( \frac{\zeta_8}{\sqrt{k}} \right)^{\dim_{\bbC} \frakh} 
		\frac{1}{\sqrt{[X : Y]}}
		\\
		&\sum_{\lambda \in X/k PY \smallsetminus \calM}
		\bm{e} \left( -\frac{1}{2Pk} \abs{\lambda}^{2} \right)
		\prod_{\alpha \in \Delta_+}
		G_{p_1, \dots, p_n} \left( \zeta_k^{\sprod{\lambda, \alpha}} \right),
	\end{align}
	where 
	\begin{itemize}
		\item denote $ \bm{e}(z) := e^{2\pi\iu z} $ for a complex number $ z $,
		\item denote $ \zeta_l^{z} := \bm{e}(z/l) $ for a positive integer $l$ and a complex number $z$, 
		\item $ X $ and $ Y $ are the weight and root lattices respectively, 
		\item $ \calM :=
		\{ \lambda \in X \mid 
		\sprod{\lambda, \alpha} \in k \Z \} $.
	\end{itemize}
\end{thm}

For a simply-laced simple Lie algebra $\frakg$, we find that the right-hand side is essentially identified with the WRT invariant. 
For a precise statement, see \cref{cor:main}.



\section{Asymptotic formulas} \label{sec:asymp}


In this section, we prepare a new asymptotic formula of false theta functions which we need to prove our main theorem.


\subsection{A new asymptotic formula} \label{subsec:asymp_extend}


To begin with, we prepare the notation for asymptotic expansion by Poincar\'{e}. 

\begin{dfn}[Poincar\'{e}]
	Let $ L $ be a positive number, $ f \colon \R_{>0} \to \bbC $ be maps, $ t $ be a variable of $ \R_{>0} $, and $ (a_n)_{n \ge -L} $ be a family of complex numbers.
	Then, we write
	\[
	f(t) \sim \sum_{n \ge -L} a_n t^{n} \text{ as } t \to +0
	\]
	if for any positive number $ M $ there exist positive numbers $ K_M $ and $ \varepsilon $ such that
	\[
	\abs{ f(t) - \sum_{-L \le n \le M} a_n t^{n} }
	\le K_M \abs{ t^{M+1}}
	\]
	for any $ 0 < t < \varepsilon $.
	In this case, we call the infinite series $ \sum_{n \ge -L} a_n t^{n} $ as the \textbf{asymptotic expansion} of  $ f(t) $ as $ t \to +0 $. 
\end{dfn}

Our new asymptotic formula is the following.
It generalizes \cite[Proposition 5.6 and Lemma 5.11]{M_plumbed}.

\begin{prop} \label{prop:asymp_lim_extend}
	Let $ (L, \sprod{\cdot, \cdot}) $ be a positive definite even lattice and
	$ \Delta \subset L \smallsetminus \{ 0 \} $ be a finite set such that the cone $ C := \linspan_{\R_{\ge 0}} \sprod{\Delta} \subset L \otimes \R $ satisfies $ -\beta \notin C $ for any $ \beta \in C \smallsetminus \{ 0 \} $.
	Let $ \calS \subset \Q^\Delta $ be a finite set such that for any $ \sigma = (\sigma_\alpha)_{\alpha \in \Delta} \in \calS $ it holds 
	$ \Delta \cdot \sigma \in L^* $,
	where $ L^* := \{ y \in L \otimes \Q \mid \sprod{x, y} \in \Z \text{ for all } x \in L \} $ is the dual lattice and 
	we denote $ \Delta \cdot x := \sum_{\alpha \in \Delta} x_\alpha \alpha \in L \otimes \R $ for $ x \in \R^\Delta $.
	Let $ \chi \colon \calS \to \bbC $ be a map, $ \delta = (\delta_\alpha)_{\alpha \in \Delta} \in \Z_{\ge 0}^\Delta $ be a vector, and $ k $ be a positive integer.
	
	Under the above settings, we prepare the following notations:
	\begin{itemize}
		\item The positive definite quadratic form $ Q(x) := \sprod{x, x}/2 $ for $ x \in L \otimes \Q $.
		\item The Gauss sum
		\[
		\calG_k(L)
		:=
		\sum_{ \lambda \in L^* / kL } 
		\bm{e} \left( -\frac{1}{k} Q(\lambda) \right).
		\]
		\item The Puiseux series
		\[
		G \left( ( q_\alpha )_{\alpha \in \Delta} \right)
		:=
		\sum_{\sigma = (\sigma_\alpha)_{\alpha \in \Delta} \in \calS} \chi(\sigma)
		\prod_{\alpha \in \Delta} q_\alpha^{\sigma_\alpha} (1-q_\alpha)^{-\delta_\alpha - 1}.
		\]
		\item The meromorphic function for complex variables 
		\[
		\varphi_k ((t_\alpha)_{\alpha \in \Delta})
		:=
		\sum_{ \lambda \in L^* / kL } 
		\bm{e} \left( -\frac{1}{k} Q(\lambda) \right)
		G \left( ( \zeta_k^{\sprod{\lambda, \alpha}} e^{-t_\alpha} )_{\alpha \in \Delta} \right).
		\]
		\item The Laurant expansion
		\[
		\varphi_k ((t_\alpha)_{\alpha \in \Delta}) =: \sum_{m \in \Z^\Delta} a_m \prod_{\alpha \in \Delta} t_\alpha^{m_\alpha},
		\]
		\item The partial theta function
		\[
		\Phi(q)
		:=
		\sum_{\sigma \in \calS} \chi(\sigma)
		\sum_{m \in \Z_{\ge 0}^\Delta} B_\delta(m) q^{Q \left( \Delta \cdot (m+\sigma) \right)},
		\]
		where for a vector $ m =  (m_\alpha)_{\alpha \in \Delta} \in \Z^\Delta $, denote
		\[
		B_\delta (m)
		:=
		\prod_{\alpha \in \Delta} \binom{m_\alpha + \delta_\alpha}{\delta_\alpha}
		\]
	\end{itemize}
	Then, we have the asymptotic formula as $ t \to +0 $
	\begin{align}
		\calG_k(L) \Phi(\zeta_k e^{-t^2})
		&\sim
		\sum_{M \in \Z} t^M
		\sum_{m \in \Z^\Delta, \, \sum_{\alpha \in \Delta} m_\alpha = M} a_m
		\restrict{ \left( \prod_{\alpha \in \Delta} \frac{\partial^{m_\alpha}}{\partial x_\alpha^{m_\alpha}} \right) e^{-Q(\Delta \cdot x)} }{x=0}.
	\end{align}	
	Here, for a Schwartz function $ g \colon \R \to \bbC $, let
	\[
	\frac{d^{-1}g}{d x^{-1}}(x) := -\int_x^\infty g(x') dx'.
	\]
\end{prop}

If we put $ \Delta $ as a basis of $ L $ in this formula, then we obtain \cite[Proposition 5.6 and Lemma 5.11]{M_plumbed}.

\begin{rem} \label{rem:gauss_sum}
	We have $\calG_k(L) = \left( \sqrt{k} \zeta_8^{-1} \right)^{\rk L} \sqrt{[L^* : L]}$ by reciprocity formula of Gauss sums (\cite[Theorem 1]{DT}).
	In particular, $\calG_k(L) \neq 0$.
\end{rem}

\begin{proof}
	To begin with, we remark that $ Q(\Delta \cdot x) \to \infty $ as $\abs{x} \to \infty$ for $ x \in \R_{\ge 0}^\Delta $ since the cone $ C $ generated by $ \Delta $ does not have lines.
	Thus, the function $ e^{-Q(\Delta \cdot x)} $ for $ x \in \R_{\ge 0}^\Delta $ is Schwartz and $ \Phi(q) $ converges.
	
	Since we have
	\begin{align}
		G \left( ( q_\alpha )_{\alpha \in \Delta} \right)
		=
		\sum_{\sigma \in \calS} \chi(\sigma)
		\sum_{m \in \Z_{\ge 0}^\Delta}
		B_\delta (m)
		\prod_{\alpha \in \Delta} q_\alpha^{m_\alpha + \sigma_\alpha}
	\end{align}
	by the binomial theorem, it holds
	\begin{align}
		\varphi_k (t)
		= \,
		&\sum_{ \lambda \in L^* / kL } 
		\bm{e} \left( -\frac{1}{k} Q(\lambda) \right)
		\sum_{\sigma \in \calS} \chi(\sigma)
		\\
		&\sum_{m \in \Z_{\ge 0}^\Delta}
		B_\delta (m)
		\prod_{\alpha \in \Delta}
		\bm{e} \left( \frac{m_\alpha + \sigma_\alpha}{k} \sprod{\lambda, \alpha} \right)
		e^{- \sum_{\alpha \in \Delta} t_\alpha (m_\alpha + \sigma_\alpha) }
		\\
		= \,
		&\sum_{\sigma \in \calS} \chi(\sigma)
		\sum_{m \in \Z_{\ge 0}^\Delta}
		B_\delta (m)
		e^{- \sum_{\alpha \in \Delta} t_\alpha (m_\alpha + \sigma_\alpha) }
		\sum_{ \lambda \in L^* / kL } 
		\bm{e} \left( 
		-\frac{1}{k} Q(\lambda) + \frac{1}{k} \sprod{\lambda, \Delta \cdot (m+\sigma)} 
		\right).
	\end{align}
	Since $ \Delta \cdot (m+\sigma) \in L^* $, we can replace $ \lambda $ by $ \lambda + \Delta \cdot (m+\sigma) $.
	Thus, we obtain
	\begin{align}
		\varphi_k (t)
		= \,
		&\sum_{\sigma \in \calS} \chi(\sigma)
		\sum_{m \in \Z_{\ge 0}^\Delta}
		B_\delta (m)
		e^{-\sum_{\alpha \in \Delta} t_\alpha (m_\alpha + \sigma_\alpha) }
		\sum_{ \lambda \in L^* / kL } 
		\bm{e} \left( 
		-\frac{1}{k} Q(\lambda) + \frac{1}{k} Q \left( \Delta \cdot (m+\sigma) \right) 
		\right)
		\\
		= \,
		&\calG_k(L)
		\sum_{\sigma \in \calS} \chi(\sigma)
		\sum_{m \in \Z_{\ge 0}^\Delta}
		B_\delta (m)
		\bm{e} \left( 
		\frac{1}{k} Q \left( \Delta \cdot (m+\sigma) \right)
		\right)
		e^{-\sum_{\alpha \in \Delta} t_\alpha (m_\alpha + \sigma_\alpha) }.
	\end{align}
	By \cite[Proposition 5.6]{M_plumbed}, we obtain the claim. 
\end{proof}


\subsection{Vanishings of asymptotic coefficients} \label{subsec:vanishing_extend}


In this subsection, we give a sufficient condition for which $ \varphi_{k} ((t_\alpha)_{\alpha \in \Delta}) $ is holomorphic at $ (t_\alpha)_{\alpha \in \Delta} = 0 $ and its constant term is ``a Gauss sum with rational function avoiding poles,'' like the expression of WRT invariants in \cref{eq:Marino}.

\begin{prop} \label{prop:value_at_0}
	Let $ (L, \sprod{\cdot, \cdot}) $ be a positive definite even lattice.
	Let $ \Delta \subset L \smallsetminus \{ 0 \} $ be a finite set, $ k $ and $ n $ be positive integers and $ p_1, \dots, p_n $ be pairwise coprime positive integers.
	Define $ P := p_1 \cdots p_n $,
	\[
	G_{p_1, \dots, p_n}(q)
	:=
	(q^{1/2} - q^{-1/2})^{2-n}
	\prod_{1 \le i \le n} (q^{1/2p_i} - q^{-1/2p_i})
	\]
	and
	\[
	\varphi_{k} ((t_\alpha)_{\alpha \in \Delta})
	:=
	\sum_{ \lambda \in L^* / kPL } 
	\bm{e} \left( -\frac{1}{kP} Q(\lambda) \right)
	\prod_{\alpha \in \Delta} G_{p_1, \dots, p_n} \left( \zeta_k^{\sprod{\lambda, \alpha}} e^{-t_\alpha} \right),
	\]
	where $ L^* := \{ y \in L \otimes \Q \mid \sprod{x, y} \in \Z \text{ for all } x \in L \} $ is the dual lattice and $ Q(x) := \sprod{x, x}/2 $ is the positive definite quadratic form for $ x \in L \otimes \Q $.
	Then, $ \varphi_{k} ((t_\alpha)_{\alpha \in \Delta}) $ is holomorphic at $ t_\alpha = 0 $ for any $ \alpha \in \Delta $ and it holds
	\begin{align}
		\varphi_{k} (0)
		=
		\sum_{ \lambda \in \left( L^* \smallsetminus \bigcup_{\alpha \in \Delta} L^*_{\alpha/k} \right) / kPL } 
		\bm{e} \left( -\frac{1}{kP} Q(\lambda) \right)
		\prod_{\alpha \in \Delta} G_{p_1, \dots, p_n} \left( \zeta_k^{\sprod{\lambda, \alpha}} \right),
	\end{align}
	where
	\[
	L^*_{\alpha/k} := \{ x \in L^* \mid \sprod{x, \alpha/k} \in \Z \}
	= \{ x \in L^* \mid \sprod{x, \alpha} \in k\Z \}.
	\]	
\end{prop}


\cref{prop:value_at_0} corresponds to \cite[Proposition 6.1]{M_plumbed}, which is proved by the same way in \cite[Theorem 4.1]{BMM_high_depth} like \cite[Proposition 4.2]{MM}.
However, its arguments used complicated elementary number theoretic methods.
In this article, we give a simple proof by focusing on the symmetry of the meromorphic function $ G \left( \{ q_\alpha \}_{\alpha \in \Delta} \right) $ like the proof of \cite[Lemma 6]{FIMT}.

\begin{proof}[Proof of $ \cref{prop:value_at_0} $]
	For a subset $ \calA \subset \Delta $, denote
	\begin{align}
		L^*_{\calA/k}
		&:=
		\bigcap_{\alpha \in \calA} L^*_{\alpha/k},
		\\
		\varphi_{k, \calA} ((t_\alpha)_{\alpha \in \Delta})
		&:=
		\sum_{ \lambda \in L^*_{\calA/k} / kPL } 
		\bm{e} \left( -\frac{1}{kP} Q(\lambda) \right)
		\prod_{\alpha \in \Delta} G_{p_1, \dots, p_n} \left( \zeta_k^{\sprod{\lambda, \alpha}} e^{-t_\alpha} \right).
	\end{align}
	By the inclusion-exclusion principle, we obtain
	\begin{align}
		&\varphi_{k} ((t_\alpha)_{\alpha \in \Delta}) + \sum_{ \emptyset \neq \calA \subset \Delta } (-1)^{\abs{\calA}} \varphi_{k, \calA} ((t_\alpha)_{\alpha \in \Delta})
		\\
		= \,
		&\sum_{ \lambda \in \left( L^* \smallsetminus \bigcup_{\alpha \in \Delta} L^*_{\alpha/k} \right) / kPL } 
		\bm{e} \left( -\frac{1}{kP} Q(\lambda) \right)
		\prod_{\alpha \in \Delta} G_{p_1, \dots, p_n} \left( \zeta_k^{\sprod{\lambda, \alpha}} e^{-t_\alpha} \right).
	\end{align}
	Here, the right-hand side is holomorphic at $ t_\alpha = 0 $ for any $ \alpha \in \Delta $ and its value at $ (t_\alpha)_{\alpha \in \Delta} = 0 $ coincides with the right hand side in the claim.
	Thus, it suffices to show that for any subset $ \emptyset \neq \calA \subset \Delta $, $ \varphi_{k, \calA} (t) $ is holomorphic at $ t_\alpha = 0 $ for any $ \alpha \in \Delta $ and has zero at $ t_\alpha = 0 $ for any $ \alpha \in \calA $. 
	We can expand
	\[
	G_{p_1, \dots, p_n}(q)
	=
	(q^{1/2} - q^{-1/2})^{2-n}
	\sum_{ (\veps_i) \in \{ \pm 1 \}^n }
	\veps_1 \cdots \veps_n
	q^{\veps_1/2p_1 + \cdots + \veps_n/2p_n}.
	\]
	For any $ \lambda \in L^*_{\calA/k} $ and $ \alpha \in \calA $, since $ \sprod{\lambda, \alpha} \in k\Z $, we have
	\begin{align}
		G_{p_1, \dots, p_n} \left( \zeta_k^{\sprod{\lambda, \alpha}} e^{-t_\alpha} \right)
		=
		(e^{-t_\alpha /2} - e^{t_\alpha /2})^{2-n}
		\sum_{ (\veps_{i}) \in \{ \pm 1 \}^n }
		&\veps_{\alpha, 1} \cdots \veps_{\alpha, n}
		e^{-t_\alpha (\veps_{1}/2p_1 + \cdots + \veps_{n}/2p_n)}
		\\
		&\bm{e} \left( 
		\frac{1}{2k} \left( -n + \frac{\veps_{1}}{p_1} + \cdots + \frac{\veps_{n}}{p_n} \right) \sprod{\lambda, \alpha}
		\right).
	\end{align}
	Thus, we obtain
	\begin{align}
		\varphi_{k, \calA} ((t_\alpha)_{\alpha \in \Delta})
		= \,
		&\sum_{\alpha \in \calA, \, (\veps_{\alpha, i}) \in \{ \pm 1 \}^n }
		\left(
		\prod_{\alpha \in \calA} (e^{-t_\alpha/2} - e^{t_\alpha/2})^{2-n}
		\prod_{1 \le i \le n } \veps_{\alpha, i} e^{-t_\alpha \veps_{\alpha, i}/2p_i}
		\right)	
		\\
		&\sum_{ \lambda \in L^*_{\calA/k} / kPL } 
		\bm{e} \left( 
		-\frac{1}{kP} Q(\lambda) 
		+ \frac{1}{2k} \sum_{\alpha \in \calA, \, 1 \le i \le n }
		\left( -1 + \frac{\veps_{\alpha, i}}{p_i} \right) \sprod{\lambda, \alpha}
		\right)
		\prod_{\beta \in \Delta \smallsetminus \calA} 
		G_{p_1, \dots, p_n} \left( \zeta_k^{\sprod{\lambda, \beta}} e^{-t_\beta} \right).
	\end{align}
	Let $ D := [L^* : L] $.
	Since $ L \supset DL^* \supset DL^*_{\calA/k} $, we have
	\begin{equation} \label{eq:varphi_calc}
		\begin{aligned}
			\varphi_{k, \calA} ((t_\alpha)_{\alpha \in \Delta})
			= \,
			&\frac{1}{[L^* : DL^*_{\calA/k}]}
			\sum_{\alpha \in \calA, \, (\veps_{\alpha, i}) \in \{ \pm 1 \}^n }
			\left(
			\prod_{\alpha \in \calA} (e^{-t_\alpha/2} - e^{t_\alpha/2})^{2-n}
			\prod_{1 \le i \le n } \veps_{\alpha, i} e^{-t_\alpha \veps_{\alpha, i}/2p_i}
			\right)
			\\
			&\sum_{ \lambda \in L^*_{\calA/k} / kPD L^*_{\calA/k} } 
			\bm{e} \left( 
			-\frac{1}{kP} Q(\lambda) 
			+ \frac{1}{2k} \sum_{\alpha \in \calA, \, 1 \le i \le n }
			\left( -1 + \frac{\veps_{\alpha, i}}{p_i} \right) \sprod{\lambda, \alpha}
			\right)
			\\
			&\prod_{\beta \in \Delta \smallsetminus \calA} 
			G_{p_1, \dots, p_n} \left( \zeta_k^{\sprod{\lambda, \beta}} e^{-t_\beta} \right).
		\end{aligned}
	\end{equation}
	Fix pairwise coprime positive integers $ \widetilde{p}_1, \dots, \widetilde{p}_n $ and integers $ e_1, \dots, e_n $ such that
	$ p_i \mid \widetilde{p}_i $, $ \widetilde{p}_1 \cdots \widetilde{p}_n = kPD $ and
	\[
	e_i \equiv
	\begin{cases}
		1 \bmod \widetilde{p}_i & \\
		0 \bmod \widetilde{p}_j & \text{ for all } j \neq i.
	\end{cases}
	\]
	By the Chinese remainder theorem, we obtain the isomorphism
	\[
	\begin{array}{ccc}
		L^*_{\calA/k} /\widetilde{p}_1 L^*_{\calA/k} \times \cdots L^*_{\calA/k} /\widetilde{p}_n L^*_{\calA/k} & \xlongrightarrow{\sim} & L^*_{\calA/k} /kPD L^*_{\calA/k} \\
		(\lambda_1, \dots, \lambda_n) & \longmapsto & e_1 \lambda_1 + \cdots + e_n \lambda_n.
	\end{array}
	\]	
	For $ \lambda_1, \dots, \lambda_n \in L^*_{\calA/k} $, we have
	\begin{align}
		\frac{1}{kP} Q \left( e_1 \lambda_1 + \cdots + e_n \lambda_n \right)
		&=
		\sum_{1 \le i \le n} \left( \frac{e_i^2}{kP} Q(\lambda_i) \right)
		+ \sum_{1 \le i < j \le n} \frac{e_i e_j}{kP} \sprod{\lambda_i, \lambda_j}
		\\
		&\equiv
		\sum_{1 \le i \le n} \left( \frac{e_i^2}{kP} Q(\lambda_i) \right)
		\bmod \Z
	\end{align}
	since $ \sprod{\lambda_i, \lambda_j} \in D^{-1} \Z $ and $ e_i e_j \in kPD \Z $.
	The last equation is invariant under $ \lambda_i \mapsto -\lambda_i $ for any fixed $ 1 \le i \le n $.
	For any $ \alpha \in \calA, \, 1 \le i \le n, \, \veps_{\alpha, i} \in \{ \pm 1 \} $, it holds
	\begin{align}
		&\frac{1}{2k} \left( -1 + \frac{\veps_{\alpha, i}}{p_i} \right) \sprod{e_1 \lambda_1 + \cdots + e_n \lambda_n, \alpha}
		\\
		= \,
		&\frac{\veps_{\alpha, i} e_i}{2k p_i} \sprod{\lambda_i, \alpha}
		- \frac{e_i }{2k} \sprod{\lambda_i, \alpha}
		+ \frac{1}{2k} \sum_{1 \le j \le n, \, j \neq i} \left( -1 + \frac{\veps_{\alpha, i}}{p_i} \right) e_j \sprod{\lambda_j, \alpha}.
	\end{align}
	On the right hand side, the second and third term is invariant under $ \lambda_i \mapsto -\lambda_i $ for any fixed $ 1 \le i \le n $ since it is a half-integer.
	The first term is invariant under $ (\veps_{\alpha, i}, \lambda_i) \mapsto (-\veps_{\alpha, i}, -\lambda_i) $ for any fixed $ 1 \le i \le n $.
	Thus, the sum for $ \lambda $ in the right hand side in \cref{eq:varphi_calc} is independent of $ (\veps_{\alpha, i}) $. 
	Therefore, by replacing $ \veps_{\alpha, i} $ as $ 1 $, we have
	\begin{equation} \label{eq:varphi_calc2}
		\begin{aligned}
			&\varphi_{k, \calA} ((t_\alpha)_{\alpha \in \Delta})
			\\
			= \,
			&\frac{1}{[L^* : DL^*_{\calA/k}]}
			\left( \prod_{\alpha \in \calA} G_{p_1, \dots, p_n} \left( e^{-t_\alpha} \right) \right)
			\\
			&\sum_{ \lambda \in L^*_{\calA/k} / kPD L^*_{\calA/k} } 
			\bm{e} \left( 
			-\frac{1}{kP} Q(\lambda) 
			+ \frac{1}{2k} \sum_{\alpha \in \calA, \, 1 \le i \le n }
			\left( -1 + \frac{1}{p_i} \right) \sprod{\lambda, \alpha}
			\right)
			\prod_{\beta \in \Delta \smallsetminus \calA} 
			G_{p_1, \dots, p_n} \left( \zeta_k^{\sprod{\lambda, \beta}} e^{-t_\beta} \right).
		\end{aligned}
	\end{equation}	
	Since $ \ord_{t=0} G(e^{-t}) = 1 $ by \cite[Lemma 2.9]{M_plumbed}, 
	the meromorphic function $ \varphi_{k, \calA} ((t_\alpha)_{\alpha \in \Delta}) $ has zero at $ t_\alpha = 0 $ for each $ \alpha \in \calA $.
	
	Finally, we prove $ \varphi_{k, \calA} ((t_\alpha)_{\alpha \in \Delta}) $ is holomorphic at $ t_\beta = 0 $ for any $ \beta \in \Delta \smallsetminus \calA $.
	For any $ \beta \in \Delta \smallsetminus \calA $, $ \varphi_{k, \calA \cup \{ \beta \}} ((t_\alpha)_{\alpha \in \Delta}) $ has zero at $ t_\beta = 0 $ by the above argument.
	The meromorphic function $ \varphi_{k, \calA} ((t_\alpha)_{\alpha \in \Delta}) - \varphi_{k, \calA \cup \{ \beta \}} ((t_\alpha)_{\alpha \in \Delta}) $ is expressed as the sum for
	$ \lambda \in L^*_{\calA/k} \smallsetminus L^*_{(\calA \cup \{ \beta \})/k} $ and
	$ G_{p_1, \dots, p_n} \left( \zeta_k^{\sprod{\lambda, \beta}} e^{-t_\beta} \right) $ is holomorphic at $ t_\beta = 0 $ for such $ \lambda $.
	Thus, $ \varphi_{k, \calA} ((t_\alpha)_{\alpha \in \Delta}) $ is hoomorphic at $ t_\beta = 0 $.
\end{proof}

By applying \cref{prop:asymp_lim_extend} as
\begin{itemize}
	\item $ k $ as $ k $,
	\item $ PL $ as $ L $,
	\item $ \sprod{\cdot, \cdot}/P $ as $ \sprod{\cdot, \cdot} $,
	\item $ P\Delta $ as $ \Delta $,
	\item define $ \calS \subset \Q^\Delta $ and $ \chi \colon \calS \to \bbC $ as
	\[
	\sum_{\sigma \in \calS} \chi(\sigma)
	\prod_{\alpha \in \Delta} q_\alpha^{\sigma_\alpha}
	:=
	\prod_{\alpha \in \Delta} 
	(-1)^n q_\alpha^{n/2 - 1}
	\prod_{1 \le i \le n} (q_\alpha^{1/2p_i} - q_\alpha^{-1/2p_i}),
	\]
	\item $ (n-3, \dots, n-3) $ as $ \delta \in \Z_{\ge 0}^\Delta $,
\end{itemize}
we obtain the following.

\begin{cor} \label{cor:value_at_0}
	Under the condition in \cref{prop:value_at_0}, define the partial theta function
	\[
	\Psi(q)
	:=
	\sum_{m \in \Z_{\ge m_0}^\Delta} q^{Q \left( \Delta \cdot m \right) / 4P}
	\prod_{\alpha \in \Delta} \chi_{m_\alpha},
	\]
	where
	\[
	\sum_{m = m_0}^\infty \chi_m q^{m/2P}
	:=
	G_{p_1, \dots, p_n} (q)
	=
	(q^{1/2} - q^{-1/2})^{2-n} \prod_{1 \le i \le n} (q_\alpha^{1/2p_i} - q_\alpha^{-1/2p_i}).
	\]
	Then, it holds
	\begin{align}
		\calG_k(L) \lim_{q \to \zeta_k} \Psi(q)
		&=
		\varphi_{k} (0)=
		\sum_{ \lambda \in \left( L^* \smallsetminus \bigcup_{\alpha \in \Delta} L^*_{\alpha/k} \right) / kPL } 
		\bm{e} \left( -\frac{1}{kP} Q(\lambda) \right)
		\prod_{\alpha \in \Delta} G_{p_1, \dots, p_n} \left( \zeta_k^{\sprod{\lambda, \alpha}} \right).
	\end{align}
\end{cor}


\section{Proof and examples of the main theorem} \label{sec:proof_ex}





Our main theorem follows from the results in \cref{sec:asymp}.

\begin{proof}[Proof of $ \cref{thm:main} $.]
	We obtain \cref{thm:main} by applying \cref{cor:value_at_0} and \cref{rem:gauss_sum} as
	\begin{itemize}
		\item $ k $ as $ k $,
		\item $ Y $ as $ L $,
		\item $ \Delta_+ $ as $ \Delta $,
		\item $\abs{x}^2/2$ as $Q(x)$.
	\end{itemize}    
\end{proof}

Finally, we give examples of our main theorem.

\begin{ex}
	Let $\frakg = \mathfrak{su}(2)$.
	Then, we have $\Delta_+ = \{ \alpha \}$ and $\abs{\alpha}^2 = 1/2$.
	Thus, we obtain $\rho = \alpha/2$, $\abs{\rho}^2 = 2$ and
	\[
	\Phi(q)
	=
	q^{-2 \phi}
	\sum_{m = m_0}^\infty \chi_{m} q^{m^{2} / 4P}.
	\]
	By \cref{rem:chi_m}, we have
	\[
	\Phi(q)
	=
	(-1)^n q^{-3\phi}
	\sum_{\veps \in \{ \pm 1 \}} \veps_1 \cdots \veps_n 
	\sum_{m' = 0}^\infty \binom{m'+n-3}{n-3}
	q^{P(\veps_1/2p_1 + \cdots \veps_1/2p_1 + n/2 - 1 + m')^2}.
	\]
	This is the same function up to simple factors as the WRT function defined by Fuji--Iwaki--Murakami--
	Terashima~\cite[Definition 1]{FIMT}.
	Andersen--Misteg\r{a}rd~\cite[Theorem 3]{Andersen-Mistegard} proved that the WRT function is equal up to simple factors to the Gukov--Pei--Putrov--Vafa invariant $\widehat{Z}(q)$ defined in \cite{GPPV}.
	In this case, we have $\Phi(q) = \widehat{Z}(q)$. 
\end{ex}

\begin{ex}
	Let $\frakg = \mathfrak{su}(3)$.
	Then, we have $\Delta_+ = \{ \alpha_1, \alpha_2, \alpha_3 \}$ and 
	\[
	S := (\sprod{\alpha_i, \alpha_j})_{1 \le i, j \le 3}
	= \pmat{2 & -1 & 1 \\ -1 & 2 & 1 \\ 1 & 1 & 2}.
	\]
	Thus, we obtain $\rho = (\alpha_1 + \alpha_2 + \alpha_3)/2$, $\abs{\rho}^2 = 2$ and
	\[
	\Phi(q)
	=
	q^{-3\phi}
	\sum_{m \in \Z_{\ge m_0}^3} \chi_{m_1} \chi_{m_2} \chi_{m_3} q^{{}^t\!m S m / 8P}.
	\]
	When $n=3$, This function is the same up to simple factors as the function $Z_{\mathrm{SU}(N)}(M)$ for $N=3$ defined by Chung~\cite[Equation 3.15]{Chung} and considered by Bringmann--Kaszian--Milas--Nazaroglu~\cite[Section 5]{BKMN_False_modular}.
	Here we remark that $S$ has eigenvalues $3$, $3$ and $0$, and thus the quadratic form ${}^t\!m S m$ is degenerate.
\end{ex}

\begin{ex}
	Let $\frakg = \mathfrak{su}(N)$ for $N \in \Z_{\ge 2}$.
	Then, we have $\Delta_+ = \{ \alpha_{i, j} \mid 1 \le i < j \le N \}$, where $\alpha_{i, j} := (0, \dots, 0, \overset{i}{1}, 0, \dots, 0, \overset{j}{-1}, 0, \dots, 0) \in \R^N$.
	Thus, we obtain
	\[
	\sprod{\alpha_{i, j}, \alpha_{i', j'}}
	=
	\begin{cases}
		2 & \text{ if } i=i', j=j', \\
		1 & \text{ if } i=i', j \neq j' \text{ or } i \neq i', j=j', \\
		-1 & \text{ if } i = j' \text{ or } i' = j, \\
		0 & \text{ otherwise}
	\end{cases}
	\]
	and
	\begin{align}
		4 \abs{\rho}^2
		&=
		\sum_{1 \le i < j \le N} \sprod{\alpha_{i, j}, \alpha_{i, j}}
		+ 2\sum_{1 \le i < j < j' \le N} \sprod{\alpha_{i, j}, \alpha_{i, j'}}
		+ 2\sum_{1 \le i < i' < j \le N} \sprod{\alpha_{i, j}, \alpha_{i', j}}
		+ 2\sum_{1 \le i < j < j' \le N} \sprod{\alpha_{i, j}, \alpha_{j, j'}} \\
		&=
		2 \# \{ (i, j) \mid 1 \le i < j \le N \}
		+ 2 \# \{ (i, j, k) \mid 1 \le i < j < k \le N \} \\
		&=
		\frac{1}{3} N(N-1)(N+1).
	\end{align}
	Let $S := (\sprod{\alpha_{i, j}, \alpha_{i', j'}})_{\substack{
			1 \le i < j \le N, \\
			1 \le i' < j' \le N
	}}$.
	Then, we obtain 
	\[
	\Phi(q)
	=
	q^{-N^2(N-1)(N+1)\phi/24}
	\sum_{m \in \Z_{\ge m_0}^{N(N-1)/2}} \chi_{m_1} \cdots \chi_{m_{N(N-1)/2}} q^{{}^t\!m S m / 8P}.
	\]
	This function is the same up to simple factors as the function $Z_{\mathrm{SU}(N)}(M)$ defined by Chung~\cite[Equation 3.15]{Chung}.
\end{ex}


\section{A relation to WRT invariants} \label{sec:WRT}


WRT invariants can be expressed as radial limits of our homological blocks as below.

\begin{cor} \label{cor:main}
	Let $ k $ be a positive integer and $ \tau_{k}^{\frakg}(M) \in \bbC $ be the WRT invariant of $ M $ with $ \frakg $.
	Then, when $\frakg$ is simply-laced, we have
	\begin{align}
		\tau_{k}^{\frakg}(M) 
		=
		&\frac{ (-1)^{\abs{\Delta_+}} \sqrt{[X \colon Y]} \zeta_8^{\dim_{\bbC} \frakg}}{\abs{W}}
		\left( \sqrt{k} \zeta_8^{-1} \right)^{\dim_{\bbC} \frakh}
		\lim_{q \to \zeta_k} \Phi(q),
	\end{align}
	where $ W $ is the Weyl group.
\end{cor}

This is obtained by our main result (\cref{thm:main}) and the following explicit formula of WRT invariants given by Mari\~{n}o~\cite[Equation (4.11) and (4.19)]{Marino} in the case when $\fg$ is simply-laced:
\begin{equation} \label{eq:Marino}
	\begin{aligned}
		\tau_{k}^{\frakg}(M) 
		=
		&\frac{ (-1)^{\abs{\Delta_+}} [X \colon Y] \zeta_8^{\dim_{\bbC} \frakg}}{\abs{W}}
		\bm{e} \left( 
		-\frac{(\dim_{\bbC} \frakg) \phi \abs{\rho}^{2} }{2k}
		\right)
		\\
		&\sum_{\lambda \in X/k PY \smallsetminus \calM}
		\bm{e} \left( -\frac{1}{2Pk} \abs{\lambda}^{2} \right)
		\prod_{\alpha \in \Delta_+}
		G_{p_1, \dots, p_n} \left( \zeta_k^{\sprod{\lambda, \alpha}} \right).
	\end{aligned}
\end{equation}

Here, we recall that we defined 
$ \calM :=
\{ \lambda \in X \mid 
\sprod{\lambda, \alpha} \in k \Z
\text{ for some } \alpha \in \Delta
\} $.
In fact, Mari\~{n}o used the set $\calM'$ given by the wall of the fundamental chamber $F_k$ with its Weyl reflections and translations by $k$ times roots. The set $\calM$ is identified with the set $\calM'$ as follows: an element $\lambda$ in $\calM'$ is written as 
\[
\lambda'=w(\omega )+k y,
\]
where $w$ is an element in the Weyl group and $\omega$ is an element in the wall of $F_k$ and $y$ is an element in the root lattice $Y$. The element $\omega$ satisfies that there exists a root $\alpha_\omega$ in $\Delta$ with $\sprod{\omega, \alpha_\omega} \in k \Z$. So, we have 
\begin{equation}
	\begin{aligned}
		\sprod{\lambda', w(\alpha_\omega)}
		&=\sprod{w(\omega), w(\alpha_\omega)}+\sprod{k y, w(\alpha_\omega)}\\
		&=\sprod{\omega, \alpha_\omega}+k \sprod{y, w(\alpha_\omega)} \in k \Z.
	\end{aligned}
\end{equation}
Then we have $\lambda' \in \calM$. Conversely, an element in $\calM$ is in the wall of some affine Weyl chamber with level $k$. Then any element $\lambda$ in $\calM$ is written as 
\[
\lambda'=w(\omega )+k y,
\]
where $w$ is an element in the Weyl group and $\omega$ is an element in the wall of $F_k$ and $y$ is an element in the root lattice $Y$, 
because each affine Weyl chamber with level $k$ is given by some Weyl reflection and translation by some root in $k Y$ of the fundamental chamber $F_k$.  


\bibliographystyle{alpha}
\bibliography{WRT_Lie_Seifert}

\end{document}